[1994/12/01]% LaTeX date must December 1994 or later
\documentclass{aomart}
\chardef\bslash=`\\ % p. 424, TeXbook
\newtheorem[{}\it]{thm}{Theorem}[section]

\newtheorem{lem}[thm]{Lemma}

\theoremstyle{definition}
\newtheorem{definition}{Definition}[section]

\newtheorem*[{}\it]{notation}{Notation}

\usepackage[all]{xy}

\newcommand{\eval}[2][\right]{\relax
  \ifx#1\right\relax \left.\fi#2#1\rvert}

\title[Relative K\"ahler-Einstein metric]{Relative K\"ahler-Einstein metric on K\"ahler varieties of positive Kodaira dimension}
\author{Hassan Jolany}
\address{Lille1 University\\
 Lille, France}
\email{hassan.jolany@math.univ-lille1.fr}
\givenname{Hassan}
\surname{Jolany}
\copyrightnote{\textcopyright~Hassan Jolany}

\subject{primary}{matsc2000}{1AB5}
\subject{primary}{matsc2000}{2FD5}
\subject{secondary}{matsc2000}{FFFF}
\subject{secondary}{matsc2000}{G25}

\begin{document}

\begin{abstract}

For projective varieties with definite first Chern class we have one type of canonical metric which is called K\"ahler-Einstein metric. But for varieties with an intermidiate Kodaira dimension we can have several different types of canonical metrics.
In this paper we introduce a new notion of canonical metric for varieties with an intermidiate Kodaira dimension. We highlight that to get $C^\infty$-solution of CMA equation of relative K\"ahler Einstein metric we need Song-Tian-Tsuji measure (which has minimal singularities with respect to other relative volume forms) be $C^\infty$-smooth and special fiber has canonical singularities. Moreover, we conjecture that if we have relative K\"ahler-Einstein metric then our family is stable in the sense of Alexeev,and  Kollar-Shepherd-Barron. By inspiring the work of Greene-Shapere-Vafa-Yau semi-Ricci flat metric, we introduce fiberwise Calabi-Yau foliation which relies in context of generalized notion of foliation. In final, we give Bogomolov-Miyaoka-Yau inequality for minimal varieties with intermediate Kodaira dimensions which admits relative K\"ahler-Einstein metric.
 
\end{abstract}

\maketitle
\tableofcontents

\section{Introduction}

Let $X_0$ be a projective variety with canonical line bundle $K\to X_0$ of Kodaira dimension $$\kappa(X_0)=\limsup\frac{\log \dim H^0(X_0,K^{\otimes \ell})}{\log\ell}$$ This can be shown to coincide with the maximal complex dimension of the image of $X_0$ under pluri-canonical maps to complex projective space, so that $\kappa(X_0)\in\{-\infty,0,1,...,m\}$.

\textbf{Lelong number}:
Let
$W\subset \mathbb C^n$
be a domain, and $\Theta$ a positive current of degree $(q,q)$ on
$W$. For a point $p\in W$
one defines
$$\mathfrak v(\Theta,p,r)=\frac{1}{r^{2(n-q)}}\int_{|z-p|<r}\Theta(z)\wedge (dd^c|z|^2)^{n-q}$$
The
Lelong number
of $\Theta$ at
$p$
is defined as

$$\mathfrak v(\Theta,p)=\lim_{r \to 0}\mathfrak v(\Theta,p,r)$$

Let $\Theta$ be the curvature of singular hermitian metric $h=e^{-u}$, one has

$$\mathfrak v(\Theta,p)=\sup\{\lambda\geq 0: u\leq \lambda\log(|z-p|^2)+O(1)\}$$

\;

Christophe Mourougane and Shigeharu Takayama, introduced the notion of relative K\"ahler metric as follows \cite{20}.

\begin{definition} Let $\pi:X\to Y$ be a holomorphic map of complex manifolds. A real d-closed $(1,1)$-form $\omega$
on $X$ is said to be a relative K\"ahler form
for $\pi$, if for every point $y\in Y$
, there exists an
open neighbourhood
$W$
of
$y$
and a smooth plurisubharmonic function
$\Psi$
on
$W$ such that $\omega+\pi^*(\sqrt[]{-1}\partial\bar\partial\Psi)$ is a K\"ahler form on $\pi^{-1}(W)$. A morphism
$\pi$ is said to be K\"ahler, if there exists a relative K\"ahler form for $\pi$, and $\pi:X\to Y$
is said to be a K\"ahler fiber space, if $\pi$ is proper, K\"ahler, and surjective with connected fibers.
\end{definition}

We consider an effective holomorphic family of complex manifolds. This means we have a holomorphic map $\pi:X\to Y$
between complex manifolds such that

1.The rank of the Jacobian of $\pi$ is equal to the dimension of $Y$ everywhere.

2.The fiber $X_t=\pi^{-1}(t)$ is connected for each $t\in Y$

3.$X_t$ is not biholomorphic to $X_{t'}$ for distinct points
$t;t'\in B$.

It is worth to mention that Kodaira showed that all fibers are dieomorphic to each other.

The relative K\"ahler form is denoted by

$$\omega_{X/Y}=\sqrt[]{-1}g_{\alpha,\bar\beta}(z,s)dz^\alpha\wedge d\bar z^\beta$$ Moreover take $\omega_X=\sqrt[]{-1}\partial\bar\partial\log \det g_{\alpha,\bar\beta}(z,y)$ on the total space $X$. The fact is $\omega_X$ in general is not K\"ahler on total space and $\omega_X|_{X_y}=\omega_{X_y}$. More precisely $\omega_X=\omega_F+\omega_H$ where $\omega_F$ is a form along fiber direction and $\omega_H$ is a form along horizontal direction. $\omega_H$
may not be K\"ahler metric in general, but $\omega_F$ is K\"ahler metric.
Now let $\omega$ be a relative K\"ahler form on $X$ and $m:=\dim X-\dim Y$, We define the relative Ricci form $Ric_{X/Y,\omega}$ of $\omega$ by

$$Ric_{X/Y,\omega}=-\sqrt[]{-1}\partial\bar\partial\log (\omega^m\wedge \pi^*|dy_1\wedge dy_2\wedge...\wedge dy_k|^2)$$
where $(y_1,...,y_k)$ is a local coordinate of $Y$, where $Y$ is a curve. See \cite{35}

Let for family $\pi:\mathcal X\to Y$
$$\rho_{y_0}:T_{y_0}Y\to H^1(X,TX)=\mathcal H_{\bar\sigma}^{0,1}(TX)$$ be the Kodaira–Spencer map for the corresponding deformation of $X$ over $Y$ at the
point $y_0\in Y$ where $\mathcal X_{y_0}=X$

If $v\in T_{y_0}Y$ is a tangent vector, say $v=\frac{\partial}{\partial y}\mid_{y_0}$ and $\frac{\partial}{\partial s}+b^\alpha\frac{\partial}{\partial z^\alpha}$ is any lift
to $\mathcal X$ along $X$, then 

$$\bar\partial \left(\frac{\partial}{\partial s}+b^\alpha\frac{\partial}{\partial z^\alpha}\right)=\frac{\partial b^\alpha(z)}{\partial z^{\bar\beta}}\frac{\partial}{\partial z^\alpha}dz^{\bar\beta}$$
is a $\bar\partial$-closed form on $X$, which represents $\rho_{y_0}(\partial/\partial y)$. 

The Kodaira-Spencer map is induced as edge homomorphism by the short exact sequence
$$0
\to T_{X/Y}\to TX
\to  \pi^*T_Y\to  0$$

This short exact sequence gives a good picture to us to run the K\"ahler Ricci flow on the relative tangent bundle.

Weil-Petersson metric when fibers are Calabi-Yau manifolds can be defined as follows\cite{6}. 
\begin{definition}
Calabi-Yau manifold is a compact K\"ahler manifold with trivial canonical bundle. The local Kuranishi family of polarized Calabi-Yau manifolds
$\mathcal X\to Y$ is smooth (unobstructed) by the Bogomolov-Tian-Todorov theorem. Let each fibers is a Calabi-Yau manifold. One can assign the unique (Ricci-flat) Yau metric $g(y)$ on $X_y$. The metric $g(y)$ induces a metric on $\wedge^{0,1}(TX)$. For
$v,w\in T_y(Y)$, one then defines the
Weil-Petersson metric
on the base
$Y$
by

$$g_{WP}(v,w)=\int_X<\rho(v),\rho(w)>_{g(y)}$$
\end{definition}

\section{Fiberwise Calabi-Yau metric}

The volume of fibers $\pi^{-1}(y)=X_y$ is a homological constant independent of $y$, and we assume that it is equal to $1$. Since fibers are Calabi-Yau manifolds so $c_1(X_y)=0$, hence there is a smooth function $F_y$ such that $Ric(\omega_y)=\sqrt{-1}\partial\bar{\partial}F_y$ and $\int_{X_y}(e^{F_y}-1)\omega_y^{n-m}=0$. The function $F_y$ vary smoothly in $y$. By Yau's
theorem there is a unique Ricci-flat K\"ahler metric $\omega_{SRF,y}$ on $X_y$ cohomologous to $\omega_0$. So there is a smooth function $\rho_y$ on $\pi^{-1}(y)=X_y$
such that $\omega_0\mid_{X_y} +\sqrt{-1}\partial\bar{\partial}\rho_y=\omega_{SRF,y}$ is the unique Ricci-flat K\"ahler metric on
$X_y$. If we normalize by $\int_{X_y}\rho_y\omega_0^n\mid_{X_y}=0$  then
$\rho_y$
varies smoothly in $y$ and defines a smooth function $\rho$ on $X$ and we let

$$\omega_{SRF}|_{X_y}=\omega_0+\sqrt{-1}\partial\bar{\partial}\rho$$
which is called as Semi-Ricci Flat metric. Such Semi-Flat Calabi-Yau metrics were first constructed by Greene-Shapere-Vafa-Yau on surfaces \cite{4}. More precisely, a closed real $(1,1)$-form $\omega_{SRF}$ on open set $U\subset X\setminus S$, (where $S$ is proper analytic subvariety contains singular points of $X$) will be called semi-Ricci flat if its restriction to each fiber $X_y\cap U$ with  $y\in f(U)$ be Ricci-flat. Notice that $\omega_{SRF}$ is positive in fiber direction, but it is still open problem that such current to be semi-positive in horizontal direction. Moreover $[\omega_{SRF}]\neq [\omega_0]$ .

\;

For the log-Calabi-Yau fibration $f:(X,D)\to Y$, such that $(X_t,D_t)$ are log Calabi-Yau varieties and central pair $(X_0,D_0)$ has simple normal crossing singularities, if $(X,\omega)$ be a K\"ahler variety with Poincaré singularities then the semi-Ricci flat metric has $\omega_{SRF}|_{X_t}$ is quasi-isometric with the following model which we call it \textbf{fibrewise Poincaré singularities.}

$$\frac{\sqrt[]{-1}}{\pi}\sum_{k=1}^n\frac{dz_k\wedge d\bar {z_k}}{|z_k|^2(\log|z_k|^2)^2}+\frac{\sqrt[]{-1}}{\pi}\frac{1}{\left(\log|t|^2-\sum_{k=1}^n\log|z_k|^2\right)^2}\left(\sum_{k=1}^n\frac{dz_k}{z_k}\wedge \sum_{k=1}^n\frac{d\bar {z_k}}{\bar {z_k}}\right)$$

We can define the same \textbf{fibrewise conical singularities.} and the semi-Ricci flat metric has $\omega_{SRF}|_{X_t}$ is quasi-isometric with the following model

$$\frac{\sqrt[]{-1}}{\pi}\sum_{k=1}^n\frac{dz_k\wedge d\bar {z_k}}{|z_k|^2}+\frac{\sqrt[]{-1}}{\pi}\frac{1}{\left(\log|t|^2-\sum_{k=1}^n\log|z_k|^2\right)^2}\left(\sum_{k=1}^n\frac{dz_k}{z_k}\wedge \sum_{k=1}^n\frac{d\bar {z_k}}{\bar {z_k}}\right)$$

In fact the previous remark will tell us that the semi Ricci flat metric $\omega_{SRF}$ has pole singularities with Poincare growth.

\;

\textbf{Remark}: Note that we can always assume the central fiber has simple normal crossing singularities(when dimension of base is one) up to birational modification and base change due to semi-stable reduction of Grothendieck, Kempf, Knudsen, Mumford and Saint-Donat as follows. 

\textbf{Theorem} (Grothendieck, Kempf, Knudsen, Mumford and Saint-Donat\cite{24})Let $ k$ be an algebraically closed field of characteristic 0 (e.g.  $ k={\mathbb{C}}$). Let $ f:X\to C$ be a surjective morphism from a  $ k$-variety $ X$ to a non-singular curve $ C$ and assume there exists a closed point $ z\in C$ such that  $ f_{\vert X\setminus f^{-1}(z)}:
X\setminus f^{-1}(z)\to C\setminus\{z\}$ is smooth. Then we find a commutative diagram

$$\displaystyle \xymatrix{ X\ar[d]_f & X\times_C C'\ar[l]\ar[d] &
X'\ar[l]_-{p}\ar[dl]^{f'}\\
C & C'\ar[l]^{\pi}
}$$
with the following properties

\;

1. $ \pi:C'\to C$ is a finite map, $ C'$ is a non-singular curve and  $ \pi^{-1}(z)=\{z'\}$.

\;

2. $ p$ is projective and is an isomorphism over  $ C'\setminus
\{z'\}$.
$ X'$ is non-singular and  $ {f'}^{-1}(z')$ is a reduced divisor with simple normal crossings, i.e., we can write  $ {f'}^{-1}(z')=\sum_i E_i$ where the $ E_i$ are 1-codimensional subvarieties (i.e., locally they are defined by the vanishing of a single equation), which are smooth and, for all  $ r$, all the intersections  $ E_{i_1}\cap\ldots\cap E_{i_r}$ are smooth and have codimension $ r$.

Now if the dimension of smooth base be bigger than one, then we don't know the semi-stable reduction and instead we can use  weak Abramovich-Karu reduction or Kawamata's unipotent reduction theorem. In fact when the dimension of base is one we know from Fujino's recent result that if we allow semi-stable reduction and MMP on the family of Calabi-Yau varieties then the central fiber will be Calabi-Yau variety. But If the dimension of smooth base be bigger than one on the family of Calabi-Yau fibers , then if we apply MMP and weak Abramovich-Karu semi-stable reduction \cite{31} then the special fiber can have simple nature. But if the dimension of base be singular then we don't know about semi-stable reduction which seems is very important for finding canonical metric along Iitaka fibration.

\section{Relative K\"ahler-Einstein metric}
\begin{definition}
Let $X$ be a smooth projective variety with $\kappa(X) \geq 0$. Then for a sufficiently
large $m > 0$, the complete linear system $|m!K_X|$ gives a rational fibration with
connected fibers $f : X \dashrightarrow  Y$.
We call $f : X  \dashrightarrow Y$ the Iitaka fibration of$ X$. Iitaka fibration is unique in the sense of birational equivalence. We may assume that $f$ is a morphism and $Y$ is
smooth. For Iitaka fibration $f$ we have

1. For a general fiber $F$, $\kappa(F) = 0$ holds.

2. $\dim Y = \kappa(Y )$.
\end{definition}

Let $X$ be a K\"ahler variety with an intermediate Kodaira dimension $\kappa(X)>0$ then we have an Iitaka fibration $\pi:X\to Y=\text{Proj} R(X,K_X)=X_{can}$ such that fibers are Calabi-Yau varieties. We set $K_{X/Y}=K_X\otimes \pi^*K_Y^{-1}$ and call it the relative canonical bundle of $\pi:X\to Y$

\begin{definition} Let $X$ be a K\"ahler variety with $\kappa(X)>0$ then the \textbf{relative K\"ahler-Einstein metric} is defined as follows 
$$Ric_{X/Y}^{h_{X/Y}^\omega}(\omega)=-\Phi\omega$$
where $\Phi$ is a fiberwise constant function, $\omega$ is the relative K\"ahler form and,
$$Ric_{X/Y}^{h_{X/Y}^\omega}(\omega)=\sqrt[]{-1}\partial\bar\partial\log(\frac{\omega^n\wedge\pi^*\omega_{can}^m}{\pi^*\omega_{can}^m})$$
and $\omega_{can}$ is a canonical metric on $Y=X_{can}$.

$$Ric_{X/Y}^{h_{X/Y}^{\omega_{SRF}}}(\omega)=\sqrt[]{-1}\partial\bar\partial\log(\frac{\omega_{SRF}^n\wedge\pi^*\omega_{can}^m}{\pi^*\omega_{can}^m})=\omega_{WP}$$
here $\omega_{WP}$ is a Weil-Petersson metric\cite{6}.

Note that if $\kappa(X)=-\infty$ then along Mori fibre space $f:X\to Y$ we can define Relative K\"ahler-Einstein metric as $$Ric_{X/Y}^{h_{X/Y}^\omega}(\omega)=\Phi\omega$$ when fibers and base are K-poly-stable. Here $\Phi$ is a fiberwise constant function. See\cite{5}

\end{definition}

Note that, if $X$ be a Calabi-Yau variety and we have a holomorphic fibre space $\pi:X\to Y$, which fibres are Calabi-Yau varieties, then we have the relative Ricci flat metric $Ric_{X/Y}(\omega)=0$. This metric is the right canonical metric on the degeneration of Calabi-Yau varieties. The complete solution of this canonical metric correspond to Monge-Ampere foliation of the fiberwise Calabi-Yau foliation and fiberwise KE stability.

\;

For the existence of K\"ahler-Einstein metric when our variety is of general type, we need to the nice deformation of K\"ahler-Ricci flow and for intermidiate Kodaira dimension we need to work on relative version of K\"ahler Ricci flow. i.e

$$\frac{\partial\omega}{\partial t}=-Ric_{X/Y}(\omega)-\Phi\omega$$
take the reference metric as $\tilde{\omega_t}=e^{-t}\omega_0+(1-e^{-t})Ric(\frac{\omega_{SRF}^n\wedge\pi^*\omega_{can}^m}{\pi^*\omega_{can}^m})$ then the version of K\"ahler Ricci flow is equivalent to the following relative Monge-Ampere equation

$$\frac{\partial\phi_t}{\partial t}=\frac{(\tilde{\omega_t}+\sqrt[]{-1}\partial\bar\partial\phi_t)^n\wedge \pi^*\omega_{can}^m}{\omega_{SRF}^n\wedge\pi^*\omega_{can}^m}-\Phi\phi_t$$

Take the relative canonical volume form $\Omega_{X/Y}=\frac{\omega_{SRF}^n\wedge\pi^*\omega_{can}^m}{\pi^*\omega_{can}^m}$ and $\omega_t=\tilde{\omega_t}+\sqrt[]{-1}\partial\bar\partial\phi_t$, then 

$$\frac{\partial\omega_t}{\partial t}=\frac{\partial\tilde{\omega_t}}{\partial t}+\sqrt[]{-1}\partial\bar\partial \frac{\partial\phi_t}{\partial t}$$

By taking $\omega_{\infty}=-Ric(\Omega_{X/Y})+\sqrt[]{-1}\partial\bar\partial\Phi\phi_{\infty}$ we obtain after using estimates

$$\log \frac{\omega_\infty^n}{\Omega_{X/Y}}-\Phi\phi_{\infty}=0$$

By taking $-\sqrt[]{-1}\partial\bar\partial$ of both sides we get

$$Ric_{X/Y}(\omega_\infty)=-\Phi\omega_\infty$$ 

and $\omega_{\infty}$ has zero Lelong number\cite{45}.

Moreover, by using higher canonical bundle formula of Kawamata, Fujino-Mori, we can have another type of canonical pair $(\omega_X,\omega_Y)$ such that $$Ric(\omega_X)=-\omega_Y+\pi^*(\omega_{WP})+[\mathcal N]$$

More explicitly on pair $(X,D)$ where $D$ is a snc divisor, we can write 

$$Ric(\omega_{(X,D)})=-\omega_{Y}+\omega_{WP}^D+\sum_P(b(1-t_P^D))[\pi^*(P)]+[B^D]$$
where $B^D$ is $\mathbb Q$-divisor on $X$ such that $\pi_*\mathcal O_X([iB_+^D])=\mathcal O_B$ ($\forall i>0$). Here $s_P^D:=b(1-t_P^D)$ where $t_P^D$ is the log-canonical threshold of $\pi^*P$ with respect to $(X,D-B^D/b)$ over the generic point $\eta_P$ of $P$. i.e., 

$$t_P^D:=\max \{t\in \mathbb R\mid \left(X,D-B^D/b+t\pi^*(P)\right)\;  \text{is sub log canonical over}\; \eta_P\}$$

For holomorphic fiber space $\pi:X\to X_{can}$, to have such pair of canonical metric, we need to have canonical bundle formula when base of fibration has canonical singularities and this is still open. In fact the canonical bundle formula of Fujino-Mori work base of CY fibration is smooth.

\textbf{Remark}:Note that the log semi-Ricci flat metric $\omega_{SRF}^D$ is not continuous in general. But if the central fiber has at worst canonical singularities and the central fiber $(X_0,D_0)$ be itself as Calabi-Yau pair, then by open condition property of Kahler-Einstein metrics, semi-Ricci flat metric is smooth in an open Zariski subset.

\textbf{Remark}:So by applying the previous remark, the relative volume form $$\Omega_{(X,D)/Y}=\frac{(\omega_{SRF}^D)^n\wedge\pi^*\omega_{can}^m}{\pi^*\omega_{can}^m\mid S\mid^2}$$
is not smooth in general, where $S\in H^0(X,L_N)$ and $N$ is a divisor which come from canonical bundle formula of Fujino-Mori. Note that Song-Tian measure is invariant under birational change

\;

Now we try to extend the Relative Ricci flow to the fiberwise conical relative Ricci flow. We define the conical Relative Ricci flow on pair $\pi:(X,D)\to Y$  where $D$ is a simple normal crossing divisor as follows 

$$\frac{\partial\omega}{\partial t}=-Ric_{(X,D)/Y}(\omega)-\Phi\omega+[N]$$
where $N$ is a divisor which come from canonical bundle formula of Fujino-Mori.

Take the reference metric as $\tilde{\omega_t}=e^{-t}\omega_0+(1-e^{-t})Ric(\frac{\omega_{SRF}^n\wedge\pi^*\omega_{can}^m}{\pi^*\omega_{can}^m})$ then the conical relative K\"ahler Ricci flow is equivalent to the following relative Monge-Ampere equation

$$\frac{\partial\phi_t}{\partial t}=\log\frac{(\tilde{\omega_t}+Ric(h_N)+\sqrt[]{-1}\partial\bar\partial\phi_t)^n\wedge \pi^*\omega_{can}^m\mid S_N\mid^2}{(\omega_{SRF}^D)^n\wedge\pi^*\omega_{can}^m}-\Phi\phi_t$$ 

Now we prove the $C^0$-estimate for this relative Monge-Ampere equation due to Tian's $C^0$-estimate

By approximation our Monge-Ampere equation, we can write

$$\frac{\partial \varphi_\epsilon}{\partial t}=\log\frac{(\omega_{t,\epsilon}+\sqrt{-1}\partial\bar{\partial}\varphi_t)^m\wedge \pi^*\omega_{can}^n\left(||S||^2+\epsilon^2\right)^{(1-\beta)}}{(\omega_{SRF}^D)^m\wedge\pi^*\omega_{can}^n}-\delta\left(||S||^2+\epsilon^2\right)^\beta-\Phi\varphi_{t,\epsilon}$$

So by applying maximal principle we get an upper bound for $\varphi_{t,\epsilon}$ as follows

$$\frac{\partial }{\partial t}\sup \varphi_\epsilon\leq \sup \log\frac{\omega_{t,\epsilon}^m\wedge \pi^*(\omega_{can})^n \left(||S||^2+\epsilon^2\right)^{(1-\beta)}}{(\omega_{SRF}^D)^m\wedge\pi^*\omega_{can}^n}-\delta\left(||S||^2+\epsilon^2\right)^\beta$$

and by expanding $\omega_{t,\epsilon}^n$ 
we have a constant $C$ independent of $\epsilon$ such that the following expression is bounded if and only if the Song-Tian-Tsuji measure be bonded, so to get $C^0$ estimate we need special fiber has mild singularities in the sense of MMP

$$\frac{\omega_{t,\epsilon}^m\wedge \pi^*(\omega_{can})^n \left(||S||^2+\epsilon^2\right)^{(1-\beta)}}{(\omega_{SRF}^D)^m\wedge\pi^*\omega_{can}^n}\approx C$$

and also $\delta\to 0$ so $\delta\left(||S||^2+\epsilon^2\right)^\beta$ is too small.
So we can get a uniform upper bound for $\varphi_{\epsilon}$. By applying the same
argument for the lower bound, and using maximal principle again, we get a $C^0$ estimate for $\varphi_{\epsilon}$. Moreover if central fiber $X_0$ has canonical singularities then Song-Tian-Tsuji measure is continuous.

So this means that we have $C^0$-estimate for relative K\"ahler-Ricci flow if and only if the central fiber has at worst canonical singularities. Note that to get $C^\infty$-estimate we need just check that our reference metric is bounded and Song-Tian-Tsuji measure is $\mathbb C^\infty$-smooth . So it just remain to see that $\omega_{WP}$ is bounded. But when fibers are not smooth in general, Weil-Petersson metric is not bounded and Yoshikawa in Proposition 5.1 in \cite{27} showed that under the some additional condition when central fiber $X_0$ is reduced and irreducible and has only canonical singularities we have $$0\leq \omega_{WP}\leq C\frac{\sqrt{-1}\mid s\mid^{2r}ds\wedge d\bar s}{\mid s\mid^2(-\log \mid s\mid)^2}$$

\section{Fiberwise Calabi-Yau foliation}
\label{lincomp}
Note that the main difficulty of the solution of $C^\infty$ for the solution of relative K\"ahler-Einstein metric is that the null direction of fiberwise Calabi-Yau metric $\omega_{SRF}$ gives a foliation along Iitaka fibration $\pi:X\to Y$ and we call it fiberwise Calabi-Yau foliation(due to H.Tsuji) and can be defined as follows $$\mathcal F=\{\theta\in T_{X/Y}|\omega_{SRF}(\theta,\bar\theta)=0\}$$
and along log Iitaka fibration  $\pi:(X,D)\to Y$, we can define the following foliation $$\mathcal F'=\{\theta\in T_{X'/Y}|\omega_{SRF}^D(\theta,\bar\theta)=0\}$$
where $X'=X\setminus D$. In fact the method of Song-Tian works when $\omega_{SRF}>0$. More precisely, in null direction, the function $\varphi$ satisfies in the complex Monge-Ampere foliation 

$$(\omega_{SRF})^\kappa=0$$

gives rise to a foliation by $X$ by complex sub-manifolds.

A complex analytic space is a topological space such that each point has an open neighborhood homeomorphic to some zero set $V(f_1,\ldots,f_k)$ of finitely many holomorphic functions in $\mathbb{C}^n$, in a way such that the transition maps (restricted to their appropriate domains) are biholomorphic functions.

\textbf{Definition}: Let $X$ be normal variety. A foliation on $X$ is a nonzero coherent subsheaf $\mathcal F \subset T_X$
satisfying

(1) $\mathcal F$ is closed under the Lie bracket, and
\;

(2) $\mathcal F$ is saturated in $T_X$ (i.e., $T_X/\mathcal F$ is torsion free). The Condition (2) above implies that $\mathcal F$ is reflexive, i.e. $\mathcal F=\mathcal F^{**}$.

\;

The canonical class $K_\mathcal F$ of $\mathcal F$ is any Weil divisor on $X$ such that $\mathcal O_X (-K_{\mathcal F} ) \cong \det(\mathcal F)$.

\begin{definition}{}
Let $\pi: X \to Y$ be a dominant morphism of normal varieties. Suppose that $\pi$ is
equidimensional. relative canonical bundle can be defined as follows $$K_{X/Y} := K_X -\pi^* K_Y$$
Let $\mathcal F$ be the foliation on $X$ induced by $\pi$, then 
$$K_{\mathcal F}=K_{X/Y}-R(\pi)$$
where $R(\pi)=\cup_D \left((\pi)^*D-((\pi)^*D)_{red}\right)$ is the ramification divisor of $\pi$. Here $D$ runs through all prime divisors on $Y$. The canonical class $K_{\mathcal F}$ of $\mathcal F$ is any Weil divisor on $X$ such that $\mathcal O_{X}(-K_{\mathcal F})\cong \det (\mathcal F):=(\wedge^r \mathcal F)^{**}$ See \cite{41}

\end{definition}

Now take a $C^\infty$ $(1,1)$-form $\omega$ on a complex manifold $X$ of complex dimension
$n$ and let $$\text{ann}(\omega)=\{W\in TX|\omega(W,\bar V)=0, \forall V\in TX\}$$

Now we have the following lemma due to Schwarz inequality \cite{16}

\begin{lem}{}If $\omega$ is non-negative then we can write,  $$\text{ann}(\omega)=\{W\in TX|\omega(W,\bar W)=0, \forall W\in TX\}$$

Moreover, if we assume $\omega^{n-1}\neq 0$ and $\omega^n=0$ then $\text{ann}(\omega)$ is subbundle of $TX$.

\end{lem}

Furthermore, we have the following straightforward lemma which make $\text{ann}(\omega)$ to be as foliation
\begin{lem}{}If $\omega$ is non-negative, $\omega^{n-1}\neq 0$, $\omega^n=0$, and $d\omega=0$, then 

 $$\mathcal F=\text{ann}(\omega)=\{W\in TX|\omega(W,\bar W)=0, \forall W\in TX\}$$
define a foliation $\mathcal F$ on $X$ and each leaf of $\mathcal F$ being a Riemann surface
\end{lem}
\;

Now Tsuji \cite{10}\cite{42} took relative form $\omega_{X/Y}$ instead $\omega$ in previous lemma and wrote it as a foliation. In my opinion Tsuji's foliation is fail to be right foliation and we need to revise it. First of all we don't know such metric $\omega_{SRF}$ is non-negative and second we must take $W\in T_{X/Y}$ in relative tangent bundle and we don't have in general $d\omega_{SRF}=0$, In fact we know just that $d_{X/Y}\omega_{SRF}=0$. Moreover $\omega_{SRF}$ is not smooth in general and it is a $(1,1)$-current with log pole singularities. 

Hence on Calabi-Yau fibration, we can introduce the following bundle
 $$\mathcal F=\text{ann}(\omega_{SRF})=\{W\in T_{X/Y}|\omega_{SRF}(W,\bar W)=0, \forall W\in T_{X/Y}\}$$
 
 in general is the right bundle to be considered and not something Tsuji wrote in \cite{14}. It is not a foliation in general. In fact it is a foliation is fiber direction and may not be a foliation in horizontal direction, but it generalize the notion of foliation. The correct solution of it as Monge-Ampere foliation still remained as open problem. 

\;

In the fibre direction, $\mathcal F$ is a foliation and we have the following straightforward theorem due to Bedford-Kalka.\cite{40}\cite{15}\cite{4}
\;

\begin{thm}{}
 Let $\mathcal L$ be a leaf of $f_*\mathcal F$, then $\mathcal L$ is a closed complex submanifold and the leaf $\mathcal L$ can be seen as fiber on the moduli map $$\eta:\mathcal Y\to \mathcal M_{CY}^D$$ where $\mathcal M_{CY}$ is the moduli space of  calabi-Yau fibers with at worst canonical singularites and 

$$\mathcal Y=\{y\in Y_{reg}|X_y\; \; \text{has Kawamata log terminal singularities}\}$$ 
 
 \end{thm}

\section{Smoothness of fiberwise integral of Calabi-Yau volume}
\label{SKE}

Let $X$ be a closed normal analytic subspace in some open subset $U$ of $\mathbb C^N$ with an isolated singularity. Take $f:X\to \Delta$ be a degeneration of smooth Calabi-Yau manifolds, then $$s\to\int_{X_s}\Omega_s\wedge\bar\Omega_s\in C^\infty$$ if and only if the monodromy $M$ acting on the cohomology of the Milnor fibre of $f$ is the
identity and the restriction map  $j: H^n(X^*) \to H^n(F)^M$ is surjective, where $X^*=X\setminus \{0\}$ and $M$ denotes monodromy acting on $H^n(F)$ and $H^n(F)^M$ is the $M$-invariant subgroup and $F$ is the Milnor fiber at zero(see Corollary 6.2. \cite{21}). In fact the $C^\infty$-smoothness of fiberwise Calabi-Yau volume $\omega_{SRF}^\kappa\wedge \pi^*\omega_Y^m$ must correspond to such information of Daniel Barlet program.

\;

Note that to get $C^\infty$-estimate for the solution of CMA along fibration $f:X\to Y$ we need to have $C^\infty$-smooth relative volume form $\Omega_{X/Y}$. So such volume forms are not unique and in fact Song-Tian-Tsuji measure has minimal singularites. If we consider a CMA equation with the relative volume form constructed by $\int_{X_s}\Omega_s\wedge \overline{\Omega_s}$, then such fiberwise integral volume forms must be smooth and in an special case when $X$ is an analytical subspace of $\mathbb C^N$ with an isolated singularities we get the $C^\infty$-smoothness of such fiberwise integral.

\;

\;

Note that, If $X_0$ only has canonical singularities, or if $X$
is smooth and $X_0$ only has isolated ordinary quadratic singularities, then
if $\pi:X\to \mathbb C^*$ be a family of degeneration of of Calabi-Yau fibers. Then the $L^2$-metric

$$\int_{X_s}\Omega_s\wedge\bar\Omega_s$$ is continuous. See Remark 2.10. of \cite{26}. 

So this fact tells us that the relative volume form is not smooth in general and finding suitable Zariski open subset such that the relative volume form (like Song-Tian-Tsuji volume form) be smooth outside of such Zariski open subset is not easy. In fact we are facing with two different singularities , one singularity arise from fiber direction near central fiber and also we have another type of singularity in horizontal direction near central fiber. So this comment tells us that  Kołodziej's $C^0$-estimate does not work for finding canonical metric along Calabi-Yau fibration.

\section{Fiberwise K\"ahler-Einstein stability}
\label{computation}
Now we use the Wang\cite{32}, Takayama\cite{29}, and Tosatti \cite{30} result for the following definition.
\begin{definition} Let $\pi:X\to B$  be a family of K\"ahler-Einstein varieties, then we introduce the new notion of stability and call it fiberwise KE-stability, if the Weil-Petersson distance $d_{WP}(B,0)<\infty$(which is equivalent to say Song-Tian-Tsuji measure is bounded near central fiber). Note when fibers are Calabi-Yau varities, Takayama, by using Tian's K\"ahler-potential for Weil-Petersson metric for moduli space of Calabi-Yau varieties showed that Fiberwise KE-Stability is as same as when the central fiber is Calabi-Yau variety with at worst canonical singularities. So this definition work when the dimension of base is one. But if the dimension of base be bigger than one, then it is better to replace boundedness of Weil-Petersson distance with boundedness of Song-Tian-Tsuji measure which seems to be more natural to me.  We mention that the Song-Tian-Tsuji measure is bounded near origin if and only if after a finite base change the Calabi-Yau family is birational to one with central fiber a Calabi-Yau variety with at worst canonical singularities.

\end{definition}

So along canonical model $\pi:X\to X_{can}$ for mildly singular variety $X$, we have $Ric_{X/X_{can}}(\omega)=-\Phi\omega$ if and only if our family of fibers be fiberwise KE-stable

Let $\pi:(X,D)\to B$
is a holomorphic submersion onto
a compact K\"ahler manifold
$B$ with $c_1(K_B)<0$ where the fibers are log Calabi-Yau manifolds and $D$ is a simple normal crossing divisor in $X$. Let our family of fibers is fiberwise KE-stable.
Then $(X,D)$
admits a unique
twisted K\"ahler-Einstein
metric $\omega_B$
solving

$$Ric(\omega_{(X,D)})=-\omega_B+\omega_{WP}^D+(1-\beta)[N]$$
where $\omega_{WP}$
is the logarithmic Weil-Petersson form on the moduli space of log Calabi-Yau fibers and $[D]$ is the current of integration over $D$.

More precisely, we have

$$Ric(\omega_{(X,D)})=-\omega_{B}+\omega_{WP}^D+\sum_P(b(1-t_P^D))[\pi^*(P)]+[B^D]$$
where $B^D$ is $\mathbb Q$-divisor on $X$ such that $\pi_*\mathcal O_X([iB_+^D])=\mathcal O_B$ ($\forall i>0$). Here $s_P^D:=b(1-t_P^D)$ where $t_P^D$ is the log-canonical threshold of $\pi^*P$ with respect to $(X,D-B^D/b)$ over the generic point $\eta_P$ of $P$. i.e., 

$$t_P^D:=\max \{t\in \mathbb R\mid \left(X,D-B^D/b+t\pi^*(P)\right)\;  \text{is sub log canonical over}\; \eta_P\}$$

and $\omega_{can}$ has zero Lelong number.

With cone angle $2\pi\beta$, $(0< \beta <1)$ along the divisor $ D$, where $h$ is an Hermitian metric on line bundle corresponding to divisor $N$, i.e., $L_N$. This equation can be solved. Take, $\omega=\omega(t)=\omega_B+(1-\beta)Ric(h)+\sqrt{-1}\partial\bar{\partial}v$ where $\omega_B=e^{-t}\omega_0+(1-e^{-t})Ric(\frac{(\omega_{SRF}^D)^n\wedge\pi^*\omega_{can}^m}{\pi^*\omega_{can}^m})$,
by using Poincare-Lelong equation, 

$$\sqrt{-1}\partial\bar{\partial}\log |s_N|_h^2=-c_1(L_N,h)+[N]$$

we have 

\begin{align*}
Ric(\omega)&=\\
&=-\sqrt{-1}\partial\bar{\partial}\log \pi_*\Omega_{(X,D)/Y}-\sqrt{-1}\partial\bar{\partial}v-(1-\beta)c_1([N],h)+(1-\beta)\{N\}\\
\end{align*}
and

\begin{align*}
\sqrt{-1}\partial\bar{\partial}\log \pi_*\Omega_{(X,D)/Y}+\sqrt{-1}\partial\bar{\partial}v&=\\
&=\sqrt{-1}\partial\bar{\partial}\log \pi_*\Omega_{(X,D)/Y}+\omega-\omega_B-Ric(h)\\
\end{align*}
Hence, by using $$\omega_{WP}^{D}=\sqrt[]{-1}\partial\bar\partial\log(\frac{(\omega_{SRF}^D)^n\wedge\pi^*\omega_{can}^m}{\pi^*\omega_{can}^m\mid S\mid^2})$$

we get 
\begin{align*}
\sqrt{-1}\partial\bar{\partial}\log \pi_*\Omega_{(X,D)/Y}+\sqrt{-1}\partial\bar{\partial}v &=\\
&=\omega_Y-\omega_{WP}^D-(1-\beta)c_1(N)\\
\end{align*}
So, 

$$Ric(\omega_{(X,D)})=-\omega_Y+\omega_{WP}^D+(1-\beta)[N]$$ 

\section{Existence of Initial K\"ahler metric along relative K\"ahler Ricci flow}

Uniqueness of the solutions of relative Kahler Ricci flow along Iitaka fibration or $\pi:X\to X_{can}$ or along log canonical model $\pi:(X,D)\to X_{can}^D$ is highly non-trivial. In fact such canonical metric is unique up to birational transformation.

Now we show how the finite generation of canonical ring can be solved by positivity theory and Analytical Minimal Model Program via K\"ahler- Ricci flow. 
\;

Now we give a relation between the existence of Zariski Decomposition and the existence of initial K\"ahler metric along relative K\"ahler Ricci flow:

Finding an initial K\"ahler metric $\omega_0$ to run the K\"ahler Ricci flow is important. Along holomorphic fibration with Calabi-Yau fibres, finding such initial metric is a little bit mysterious. In fact, we show that how the existence of initial K\"ahler metric is related to finite generation of canonical ring along singularities.

Let $\pi:  X\to Y$ be an Iitaka fibration of projective varieties
$X,Y$,(possibily singular) then is there always the following decomposition
 $$K_Y+\frac{1}{m!}\pi_*\mathcal O_X(m!K_{X/Y})=P+N$$ 
  where $P$ is semiample and $N$ is effective.  The reason is that, If $X$ is smooth projective variety, then as we mentioned before, the canonical ring $R(X,K_X)$ is finitely generated. We may thus assume that $R(X,kK_X)$ is generated in degree 1 for some $k>0$. Passing to a log resolution of $|kK_X|$ we may assume that $|kK_X|=M+F$ where $F$ is the fixed divisor and $M$ is base point free and so $M$ defines a morphism $f:X\to Y$ which is the Iitaka fibration. Thus $M=f^*O_Y(1)$ is semiample and $F$ is effective.

In singular case, if $X$ is log terminal. By using Fujino-Mori's higher canonical bundle formula, after resolving $X'$, we get a morphism $X'\to Y'$ and a klt pair $K_Y'+B_Y'$. The Y described above is the log canonical model of $K_Y'+B_Y'$ and so in fact (assuming as above that $Y'\to Y$ is a morphism), then $K_Y'+B_Y'\sim_{\mathbb Q} P+N$ where $P$ is the pull-back of a rational multiple of $O_Y(1)$ and $N$ is effective (the stable fixed divisor). If $Y'-\to Y$ is not a morphism, then P will have a base locus corresponding to the indeterminacy locus of this map.(Thanks of Hacon answer to my Mathoverflow question \cite{22} which is due to E.Viehweg \cite{23} )

So the existence of Zariski decomposition is related to the finite generation of canonical ring (when $X$ is smooth or log terminal). Now if such Zariski decomposition exists then, there exists a singular hermitian metric $h$, with semi-positive Ricci curvature $\sqrt[]{-1}\Theta_h$ on $P$, and it is enough to take the initial metric $\omega_0=\sqrt[]{-1}\Theta_h+[N]$ or $\omega_0=\sqrt[]{-1}\Theta_h+\sqrt[]{-1}\delta\partial\bar\partial \|S_N\|^{2\beta}$ along relative K\"ahler Ricci flow
$$\frac{\partial \omega(t)}{\partial t}=-Ric_{X/Y}(\omega(t))-\Phi\omega(t)$$
with log terminal singularities.

So when $X,Y$ have at worst log terminal singularities(hence canonical ring is f.g and we have initial K\"ahler metric to run K\"ahler Ricci flow with starting metric $\omega_0$)  and central fibre is Calabi-Yau variety, and $-K_Y<0$, then all the fibres are Calabi-Yau varieties and the relative K\"ahler-Ricci flow converges to $\omega$  which satisfies in $$Ric_{X/Y}(\omega)=-\Phi\omega$$

\textbf{Remark}: The fact is that the solutions of relative K\"ahler-Einstein metric or Song-Tian metric $Ric(\omega_X)=-\omega_Y+f^*\omega_{WP}+[\mathcal N]$ may not be $C^\infty$. In fact we have $C^\infty$ of solutions if and only if the Song-Tian measure or Tian's K\"ahler potential be $C^\infty$. Now we explain that under some following algebraic condition we have $C^\infty$-solutions for $$Ric(\omega_X)=-\omega_Y+f^*\omega_{WP}+[\mathcal N]$$ along Iitaka fibration. We recall the following Kawamata's theorem \cite{17}.
\begin{thm}
Let $f : X\to B$ be a surjective morphism of smooth projective
varieties with connected fibers. Let $P=\sum_jP_j$, $Q=\sum_lQ_l$, be normal crossing
divisors on $X$ and $B$, respectively, such that $f^{-1}(Q)\subset P$ and $f$ is smooth over $B\setminus Q$. Let $D =\sum_j d_jP_j$ be a $\mathbb Q$-divisor on $X$, where $d_j$ may be positive, zero or negative, which satisfies the following conditions A,B,C:

\;

\textbf{A}) $D = D^h + D^v$ such that any irreducible component of $D^h$ is mapped surjectively
onto $B$ by $f$ , $f : Supp(D^h) \to B$ is relatively normal crossing over $B \setminus Q$, and
$f(Supp(Dv ))\subset  Q$. An irreducible component of $D^h$ (resp. $D^v$ ) is called horizontal
(resp. vertical)

\;

\textbf{B})$d_j < 1$ for all $j$

\;

\textbf{C}) The natural homomorphism $\mathcal O_B \to f_*\mathcal O_X(\lceil -D\rceil)$
is surjective at the generic point of $B$.

\textbf{D}) $K_X + D\sim_{\mathbb Q} f^*
(K_B + L)$ for some $\mathbb Q$-divisor $L$ on $B$.

Let

\begin{align*}
f^*Q_l &= \sum_jw_{lj}P_j \\
 \bar{d_j} &=\frac{d_j+w_{lj}-1}{w_{lj}},\; \text{if}\;  f(P_j)=Q_l \\
\delta_l  &= max\{ \bar{d_j}; f(P_j)=Q_l \}.\\
\Delta  &= \sum_l \delta_l Q_l .\\
M  &= L-\Delta .\\
\end{align*}

Then $M$ is nef.
\end{thm}

The following theorem is straightforward from Kawamata's theorem 
\begin{thm}Let $d_j < 1$ for all $j$ be as above in Theorem 0.11, and fibers be log Calabi-Yau pairs, then $$\int_{X_s\setminus D_s}(-1)^{n^2/2}\frac{\Omega_s\wedge\overline{\Omega_s}}{\mid S_s\mid^2}$$  is continuous on a nonempty Zariski open subset of $B$. 

\end{thm}

Since the inverse of volume gives a singular hermitian line bundle, we have the following theorem from Theorem 0.11

\begin{thm}Let $K_X + D\sim_{\mathbb Q} f^*
(K_B + L)$ for some $\mathbb Q$-divisor $L$ on $B$ and

\begin{align*}
f^*Q_l &= \sum_jw_{lj}P_j \\
 \bar{d_j} &=\frac{d_j+w_{lj}-1}{w_{lj}},\; \text{if}\;  f(P_j)=Q_l \\
\delta_l  &= max\{ \bar{d_j}; f(P_j)=Q_l \}.\\
\Delta  &= \sum_l \delta_l Q_l .\\
M  &= L-\Delta .\\
\end{align*}

Then 

$$\left(\int_{X_s\setminus D_s}(-1)^{n^2/2}\frac{\Omega_s\wedge\overline{\Omega_s}}{\mid S_s\mid^2}\right)^{-1}$$
is a continuous hermitian metric on the $\mathbb Q$-line bundle $K_B + \Delta$ when fibers are log Calabi-Yau pairs.

\end{thm}

\section{Stable family and Relative K\"ahler-Einstein metric}
\label{s:comp}

For compactification of the moduli spaces of polarized varieties Alexeev,and  Kollar-Shepherd-Barron,\cite{25} started a program by using new notion of moduli space of "stable family". They needed to use the new class of singularities, called semi-log canonical singularities. 

\;

Let $X$ be an equidimensional algebraic variety that
satisfies Serre's $S_2$ condition and is normal crossing in codimension one. Let $\Delta$ be an effective
$\mathbb R$-divisor whose support does not contain any irreducible components of the conductor of $X$.
The pair $(X, \Delta)$ is called a semi log canonical pair (an slc pair, for short) if

\;

(1) $K_X + \Delta$ is $\mathbb R$-Cartier;

\;

(2) $(X^v, \Theta)$ is log canonical, where $v: X^v \to X$ is the normalization and $K_{X^v}+\Theta =v^*(KX +\Delta)$

Note that, the conductor $\mathcal C_X$ of $X$ is the subscheme defined by, $\mathfrak{cond}_X:={Hom}_{\mathcal O_X}(v_*\mathcal O_{X^v},\mathcal O_X)$.

\;

A morphism $f : X \to B$ is called a weakly stable family if it
satisfies the following conditions:
\;

1. $f$ is flat and projective
\;

2. $\omega_{X/B}$ is a relatively ample $\mathbb Q$-line bundle
\;

3. $X_b$ has semi log canonical singularities for all $b \in B$

\;
 A weakly stable family $f : X \to B$ is called a stable family if
it satisfies Kollar’s condition, that is, for any $m \in \mathbb N$

$$\omega^{[m]}_{X/B}|_{X_b}\cong \omega_{X_b}^{[m]}
.$$
\;

Note that, if the central fiber be Gorenstein and stable variety, then all general fibers are stable varieties, i.e, stability is an open condition
\;

\textbf{Conjecture}: Weil-Petersson metric (or logarithmic Weil-Petersson metric)on stable family is semi-positive as current and such family has finite distance from zero i.e $d_{WP}(B,0)<\infty$ when central fiber is stable variety also.

\;

Moreover we predict the following conjecture holds true.

\;

\textbf{Conjecture}: Let $f: X \to B$ is a stable family of polarized Calabi-Yau varieties, and let $B$ is a smooth disc. then if the central fiber be stable variety as polarized Calabi-Yau variety, then we have following canonical metric on total space. $$Ric(\omega_X)=-\omega_B+f^*(\omega_{WP})+[N]$$

Moreover, if we have such canonical metric then our family of fibers is stable.

\;

We predict that if the base be singular with mild singularites of general type(for example $B=X_{can}$) then we have such canonical metric on the stable family

Now the following formula is cohomological characterization of Relative K\"ahler-Ricci flow due to Tian
\begin{thm} The maximal time existence $T$ for the solutions of relative K\"ahler Ricci flow is $$T=\sup\{t\mid e^{-t}[\omega_0]+(1-e^{-t})c_1(K_{X/Y}+D)\in \mathcal K((X,D)/Y)\}$$
where  $\mathcal K\left((X,D)/Y\right)$ denote the relative K\"ahler cone of $f:(X,D)\to Y$
\end{thm}

Now take we have holomorphic fibre space $f:X\to Y$ such that fibers and base are Fano K-poly stable, then we have the relative K\"ahler-Einstein metric $$Ric_{X/Y}(\omega)=\Phi\omega$$ we need to work on relative version of K\"ahler Ricci flow. i.e

$$\frac{\partial\omega}{\partial t}=-Ric_{X/Y}(\omega)+\Phi\omega$$
take the reference metric as $\tilde{\omega_t}=e^{t}\omega_0+(1+e^{t})Ric(\frac{(\omega_{SKE}^{n})\wedge\pi^*\omega_{Y}^m}{\pi^*\omega_{Y}^m})$ then the version of K\"ahler Ricci flow is equivalent to the following relative Monge-Ampere equation

$$\frac{\partial\phi_t}{\partial t}=\frac{(\tilde{\omega_t}+\sqrt[]{-1}\partial\bar\partial\phi_t)^n\wedge \pi^*\omega_{Y}^m}{(\omega_{SKE}^{n})\wedge\pi^*\omega_{Y}^m}+\Phi\phi_t$$

where $\omega_Y$ is the K\"ahler-Einstein metric corresponding to $Ric(\omega_Y)=\omega_Y$ and $\omega_{SKE}$ is the fiberwise Fano K\"ahler-Einstein metric.

In fact the relative volume form is $\Omega_{X/Y}=\frac{(\omega_{SKE}^{n})\wedge\pi^*\omega_{Y}^m}{\pi^*\omega_{Y}^m}$ and we have the following relative Monge-Ampere equation 

$$\frac{\partial\phi_t}{\partial t}=\frac{(\tilde{\omega_t}+\sqrt[]{-1}\partial\bar\partial\phi_t)^n}{\Omega_{X/Y}}+\Phi\phi_t$$

Hence from $Ric_{X/Y}(\omega)=\Phi\omega$. Moreover we must develop canonical bundle formula when fibers are K-stable and if we have such formula then we obtain $Ric(\omega_X)=\omega_Y+f^*(\omega_{WP})+[\mathcal S]$, for the Weil-Petersson metric $\omega_{WP}$ on the base which this metric is correspond to canonical metric on moduli part of family of fibers, $\omega_{WP}=\int_{X/Y}c_1(K_{X/Y},h)^{n+1}$. 

\textbf{Remark}: Note that we still don't know canonical bundle type formula along Mori-fiber space. So finding explicit Song-Tian type metric on pair $(X,D)$ along Mori fiber space when base and fibers are K-poly stable is not known yet.

\;
\;
\;

\textbf{Conjecture}:Let $\pi: X\to B$ is smooth, and every
$X_t$ is K-poly stable.
Then the plurigenera $P_m(X_t)=\dim H^0(X_t,-mK_{X_t})$ is independent of $t\in B$ for any $m$.

Idea of proof. We can apply the relative K\"ahler Ricci flow method for it. In fact if we prove that $$\frac{\partial \omega(t)}{\partial t}=-Ric_{X/Y}(\omega(t))+\Phi\omega(t)$$
has long time solution along Fano fibration such that the fibers are K-poly stable then we can get the invariance of plurigenera in the case of K-poly stability

\;

\section{Bogomolov-Miyaoka-Yau inequality for minimal varieties with intermediate Kodaira dimension}

From the differential geometric proof of Yau \cite{44} and the algebraic proof of Miyaoka \cite{43} for minimal varieties of general type $\kappa (X)=\dim X$, we know that by using K\"ahler Ricci flow method we can get the following inequality $$(-1)^nc_1^n(X)\leq (-1)^n\frac{2(n+1)}{n}  c_1^{n-2}(X)c_2(X)$$

 So we can extend this idea for the Bogomolov-Miyaoka-Yau inequality for minimal
 varieties with an intermediate Kodaira dimension  $0<\kappa(X)<\dim
 X$
 
 So, we have the following inequality as soon as relative K\"ahler Ricci flow has $C^{\infty}$-solution:

 $$\left(\frac{2(n-m+1)}{n-m}c_2(\mathcal T_{X/X_{can}})-c_1^2(\mathcal
 T_{X/X_{can}})\right).[\Phi\omega]^{n-2}\geq 0$$

 where $\omega$ is a relative K\"ahler form on the minimal projective variety
 $X=X_{min}$ and $X_{can}=\text{Proj}\bigoplus_{m\geq 0}H^0(X,K_X^{m})$
 is the canonical model of $X$ (here  $\mathcal T_{X/X_{can} 
}=Hom(\Omega^1_{X/X_{can}}, \mathcal O_X)$ mean relative tangent sheaf) via
 Iitaka fibration $\pi: X\to X_{can}$.

Certainly we must require stability in order that this inequality holds true. The stability must be equivalent with the fact that the following flow $C^\infty$-converges in $C^\infty$$$\frac{\partial\omega(t)}{\partial t}=-Ric_{X/X_{can}}(\omega(t))-\Phi\omega(t)$$

Here $Ric_{X/X_{can}}=dd^c\log \Omega_{X/X_{can}}$(where $\Omega_{X/X_{can}}$ is the relative volume form) means relative Ricci form and $\Phi$ is fiberwise constant function. Note that if such relative K\"ahler Ricci flow has solution then $K_{X/X_{can}}$ is psudo-effective
I think that the analytical minimal model program can prove this. 

\;

In fact, if we have relative K\"ahler-Einstein metric $Ric_{X/X_{can}}\omega=-\Phi\omega$ ,  then Bogomolov-Miyaoka-Yau inequality for minimal varieties with intermediate Kodaira dimension $0<\kappa (X)<\dim X$ holds true.

\section{Canonical metric on foliations}

In fact we can study the canonical metric on foliations on its canonical model of projective varieties. We have minimal model program on foliations developed by Michael McQuillan \cite{49} and we can extend Song-Tian program on foliations. But in general study of canonical metric on foliations is more complicated. For example abundance conjecture is not true on foliations. We need to the analytical surgery by using Partial K\"ahler Ricci flow,or mixed scalar curvature introduced by Vladimir Rovenski and Vladimir Sharafutdinov which must be compatible with algebraic surgery (Minimal Model program on foliations). So when for the foliation $\mathcal F$ of general type, we have $c_1(\mathcal F)<0$ the right canonical metric can be as same as the K\"ahler-Einstein metric but instead Ricci curvature we must use Partial Ricci curvature and leafwise constant to design such canonical metric. The fact is that we need new techniques to get $C^\infty$ solution and continuity method does not work!.

\section{Acknowledgements}
I would like to thanks of my Ph.D. advisor, professor Gang Tian for proposing this project to me and invited me to visit him in department of Mathematics of Princeton university.

\end{document}